\documentclass[preprint,3p,12pt,pdf]{elsarticle}

\usepackage{lineno,hyperref}
\modulolinenumbers[5]

\usepackage{hyperref}

\usepackage{epstopdf}

\usepackage{tikz}
\usetikzlibrary{arrows}

\usepackage{pst-node}
\usepackage{tikz-cd}

\usepackage[shellescape]{gmp}

\usepackage{mathrsfs}
\usepackage{amssymb}
\usepackage{amsfonts}
\usepackage{latexsym}
\usepackage{mathtools}
\usepackage{xcolor}
\usepackage{wrapfig}
\usepackage{floatflt}

\usepackage{mathtools}
\usepackage{extarrows}

\usepackage{graphicx}

\usepackage{subcaption}

\def\lb{\label}

\newcommand{\er}[1]{\textrm{(\ref{#1})}}

\newtheorem{theorem}{\bf Theorem}[section]

\def\g{\gamma}

           \def\mJ{{\mathscr J}}
        
  \def\cL{{\mathcal L}}       
 \def\cM{{\mathcal M}}       \def\mM{{\mathscr M}}

\def\r{\rho}           
\def\s{\sigma}         
         
\def\t{\tau}

\def\P{\Psi}

\def\ve{\varepsilon}   \def\vt{\vartheta}    \def\vp{\varphi}    

\def\Z{{\mathbb Z}}    \def\R{{\mathbb R}}   \def\C{{\mathbb C}}    
    \def\N{{\mathbb N}}   

\def\lt{\biggl}                  \def\rt{\biggr}
\def\ol{\overline}               \def\wt{\widetilde}

\let\ge\geqslant                 \let\le\leqslant

\def\iy{\infty}
\def\sm{\setminus}               \def\es{\emptyset}

\def\el2{\ell^{\,2}}             \def\1{1\!\!1}

\def\Im{\mathop{\mathrm{Im}}\nolimits}

\def\Re{\mathop{\mathrm{Re}}\nolimits}

\def\sign{\mathop{\mathrm{sign}}\nolimits}

\newtheorem{corollary}[theorem]{\bf Corollary}

\let\ge\geqslant
\let\le\leqslant

\newcommand{\ca}{\begin{cases}}
\newcommand{\ac}{\end{cases}}
\newcommand{\ma}{\begin{pmatrix}}
\newcommand{\am}{\end{pmatrix}}
\def\eq{\begin{equation}}
\def\qe{\end{equation}}
\def\[{\begin{equation}}
\def\]{\end{equation}}

\begin{document}

\begin{frontmatter}

\title{Exponential stochastic compression of one-dimensional space and 146 percent
}

\date{\today}

\author
{Anton A. Kutsenko}

\address{Jacobs University, 28759 Bremen, Germany; email: akucenko@gmail.com}

\begin{abstract}
Exponential stochastic compression is the process when every second cell of an infinite chain may increase its weight merging randomly with left, right, or both neighboring cells. The total mass conservation is assumed. After that, merged cells fill the empty space, compressing the chain twice. They may fill empty spaces in two different ways: (I) using shifts only, i.e. preserving the order; (II) using shifts and random permutations. Compressing the initial homogeneous chain with cell weights $1$ many times, we compute final densities $\rho_i$ of cells with weight $i=1,2,3,...$. The main result is that $\rho_i/\rho_1=i$ in the ordered case (I), and $\rho_i/\rho_1\approx1.464910...(i-1/4)$ in the disordered case (II). The multiplier in the disordered case has a fractal nature. The compression of initially inhomogeneous chains and rescaled continuous densities are also discussed.
\end{abstract}

\begin{keyword}
compression, stochastic, fractals, polynomial dynamics
\end{keyword}


\end{frontmatter}


{\section{Introduction}\lb{sec1}}

Current research is motivated by the analysis of a primitive one-step model of group formation, proposed in \cite{K1}. I extend this model by adding arbitrary large number of steps. In this case, the explicit analysis is still possible. It is based on ideas from the theory of branching processes and polynomial dynamics. However, any analytic results regarding the non-trivial dynamical system with stochastic input proposed in \cite{K1} still seem unreachable. In fact, I would not write this new work if I found the constant $1.464910...$ in the knowledge base. On the other hand, this new work contains several analytic aspects in an unusual combination worthy of exposition. And finally, these 146 percent are not related to corrupt or compromised elections in 2011.

I have tried to make the presentation easy to read, following the spirit of classic monographs in the theory of probability such as \cite{F}. A very good introduction to the polynomial and holomorphic dynamics with a nice explanation of fundamental ideas can be found in \cite{M}. Subjectively for me, work \cite{K2}, where some infinite product expansions of solutions of Poincar\'e and Shcr\"oder-type functional equations are discussed, was also very helpful in this research. Investigating conformal mappings related to the solutions of Shcr\"oder-type functional equations, I found that some of its components having a fractal nature can be well approximated by constant functions with a very small error. This idea is used in the current research to find the constant $1.464910...$. It is useful to note that near-constancy oscillations of other fractal functions are considered in, e.g., \cite{CG}, see also references therein, in a different context. In the references, where all the papers are good enough, due to my subjective taste, I would like to highlight \cite{H}, \cite{BB}, and especially \cite{DIL} for its light style.

Consider the infinite chain $\Z$, where at each of the cell the value $1$ is placed. Then each even cell becomes active. It moves independently left or right with the probability $1/2$ and its value is added to the value in the corresponding neighboring odd cell. After that, odd cells fill the empty spaces, i.e. compressing the chain twice. Thus, there are no empty cells in the chain after the compression. All the cells have non-zero weights. We consider two ways of filling empty spaces: ordered compression, where the only simple shifts without permutations are allowed; disordered compression, where odd cells fill the empty spaces randomly using permutations. 

There is another interpretation of ordered and disordered stochastic compression. Imagine an infinite chain of citizens, each of which has exactly one coin. Then every second citizen gives her/his money randomly to the left or right neighbor with equal probability and leaves the chain. The remaining citizens in the chain fill the places vacated after the departure of neighbors. They can do this by keeping order, or by moving randomly with permutations. An example of one step of compression is given in Fig. \ref{fig0}. The process of compression is repeated many times. As will be shown below, the random permutations significantly affect coin accumulation. The explicit analysis of such a seemingly simple problem involves various techniques of complex analysis related to fractal Julia sets. The ``citizens and coins" interpretation is clear enough, but we return to the original formulation in terms of the chain of cells and their weights.
\begin{figure}[h]
	\center{\includegraphics[width=0.9\linewidth]{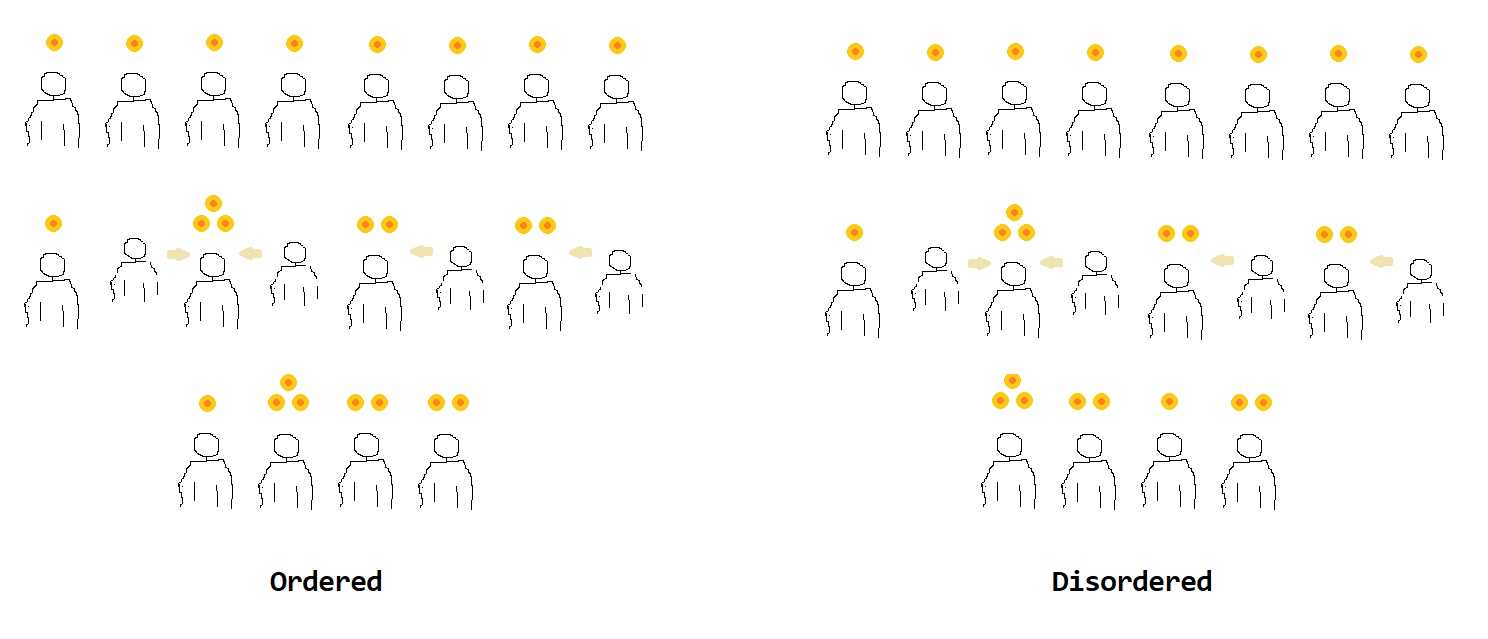}}
	\caption{Ordered and disordered stochastic compression. Interpretation with citizens and coins.}\lb{fig0}
\end{figure}

Before we proceed with analytic results, let us try numerical simulations. Analyzing results in Table \ref{Tab1}, we see that the exponential ordered stochastic compression (OSC) leads to the simple linear growth of densities of cells with different weights, whereas the exponential disordered stochastic compression (DSC) leads to a more complex, perhaps also linear, growth of densities.

\begin{table}[h]
\begin{center}
    \begin{tabular}{ | l | l | l | l | l | l | l | l | l |}
    \hline
    i & $1$ & $2$ & $3$ & $4$  & $5$ & $6$ & $7$ & $8$ \\ \hline
    $N_i$ & $6.55\cdot10^8$  & $13.1\cdot10^8$  & $19.65\cdot10^8$ & $26.2\cdot10^8$  & $19.65\cdot10^8$  & $13.1\cdot10^8$  & $6.5\cdot10^8$ & $0$ \\ \hline
    $\frac{N_i}{N_1}$ & $1.00000$ & $2.00018$ & $3.00015$ & $4.00041$  & $3.00003$  & $2.00021$  & $1.00001$ & $0$ \\ 
    \hline
    \end{tabular}
\\
\bigskip
(a) OSC, $2$ steps
\bigskip
\\
    \begin{tabular}{ | l | l | l | l | l | l | l | l | l |}
    \hline
    i & $1$ & $2$ & $3$ & $4$  & $5$ & $6$ & $7$ & $8$ \\ \hline
    $N_i$ & $2413$  & $4972$  & $7463$ & $10172$  & $12322$  & $14996$  & $17561$ & $19885$ \\ \hline
    $\frac{N_i}{N_1}$ & $1.00000$ & $2.06051$ & $3.09283$ & $4.21550$  & $5.10651$  & $6.21467$  & $7.27766$ & $8.24078$ \\ 
    \hline
    \end{tabular}
\\    
    \bigskip
    (b) OSC, $8$ steps
\\
 \bigskip 
 \begin{tabular}{ | l | l | l | l | l | l | l | l | l |}
    \hline
    i & $1$ & $2$ & $3$ & $4$  & $5$ & $6$ & $7$ & $8$ \\ \hline
    $N_i$ & $6.5\cdot10^8$  & $16.2\cdot10^8$  & $20\cdot10^8$ & $22\cdot10^8$  & $19.2\cdot10^8$  & $11.5\cdot10^8$  & $6.1\cdot10^8$ & $2.3\cdot10^8$ \\ \hline
    $\frac{N_i}{N_1}$ & $1.00000$ & $2.48012$ & $3.05856$ & $3.40276$  & $2.94165$  & $1.75454$  & $0.94146$ & $0.36279$ \\ 
    \hline
    \end{tabular}
\\
\bigskip
(c) DSC, $2$ steps
\bigskip
\\
    \begin{tabular}{ | l | l | l | l | l | l | l | l | l |}
    \hline
    i & $1$ & $2$ & $3$ & $4$  & $5$ & $6$ & $7$ & $8$ \\ \hline
    $N_i$ & $2508$  & $6571$  & $9462$ & $13741$  & $17573$  & $20134$  & $24709$ & $28157$ \\ \hline
    $\frac{N_i}{N_1}$ & $1.00000$ & $2.62002$ & $3.77273$ & $5.47887$  & $7.00678$  & $8.02791$  & $9.85207$ & $11.22687$ \\ 
    \hline
    \end{tabular}
\\    
    \bigskip
    (d) DSC, $8$ steps   
\end{center}
\caption{The quantities $N_i$ of cells with weights $i$ in a segment of homogeneous chain of length $2^{26}\cdot 5^4\approx42\cdot10^9$ after $2$ and $8$ steps of ordered and disordered stochastic compression, in numerical simulations.}\lb{Tab1}
\end{table}
Now we know more or less what we will prove. The density $\r_i$ of the cells with weight $i$ in the infinite chain $\Z$ can be naturally defined as
\[\lb{001}
 \r_i=\lim_{N\to\iy}\frac1{2N}\#\{j:\ x_j=i,\ j\in[-N,...,N-1]\},
\]
where $\#$ denotes the number of elements in the set.

\begin{theorem}\lb{T1} We take the infinite chain $\Z$, where $1$ is placed at each of the cells.

i) After $N$ steps of the exponential ordered stochastic compression, the densities $\r_i$ of the cells with weights $i$ are
\[\lb{002}
 \r_i=\frac1{4^N}\ca 
          i,& i=1,...,2^N,\\
          2^{N+1}-i,& i=2^N+1,...,2^{N+1}-1,\\
          0,& i\ge 2^{N+1}.  
      \ac
\]

ii) After $N$ steps of the exponential disordered stochastic compression, the densities $\r_i$ of the cells with weights $i$ satisfy
\[\lb{003}
 \sum_{i=1}^{\iy}\r_iz^i=\underbrace{P\circ...\circ P}_{N}(z),\ \ z\in\C,
\]
where $P(z)=\frac14(z+2z^2+z^3)$.
\end{theorem}

While the result i) of Theorem \ref{T1} looks more simple than ii), the Proof of i) is more complex. To find approximations of $\r_i$ in DSC ii), it is natural to introduce the analytic function
\[\lb{004}
 \Phi(z)=\lim_{N\to\iy}4^N\underbrace{P\circ...\circ P}_{N}(z),
\]
which satisfy the Schr\"oder-type functional equation
\[\lb{005}
 \Phi(P(z))=\frac14\Phi(z),\ \ \ \Phi'(0)=1.
\]
The Koenig's Theorem guarantees the existence of $\Phi(z)$, since $0$ is an attracting point with the multiplier $\frac14$: $P(0)=0$ and $P'(0)=\frac14$. The domain of definition of $\Phi(z)$ coincides with the filled Julia set related to the polynomial $P(z)$. The corresponding Julia set is illustrated in Fig. \ref{fig1}.
\begin{figure}[h]
	\center{\includegraphics[width=0.5\linewidth]{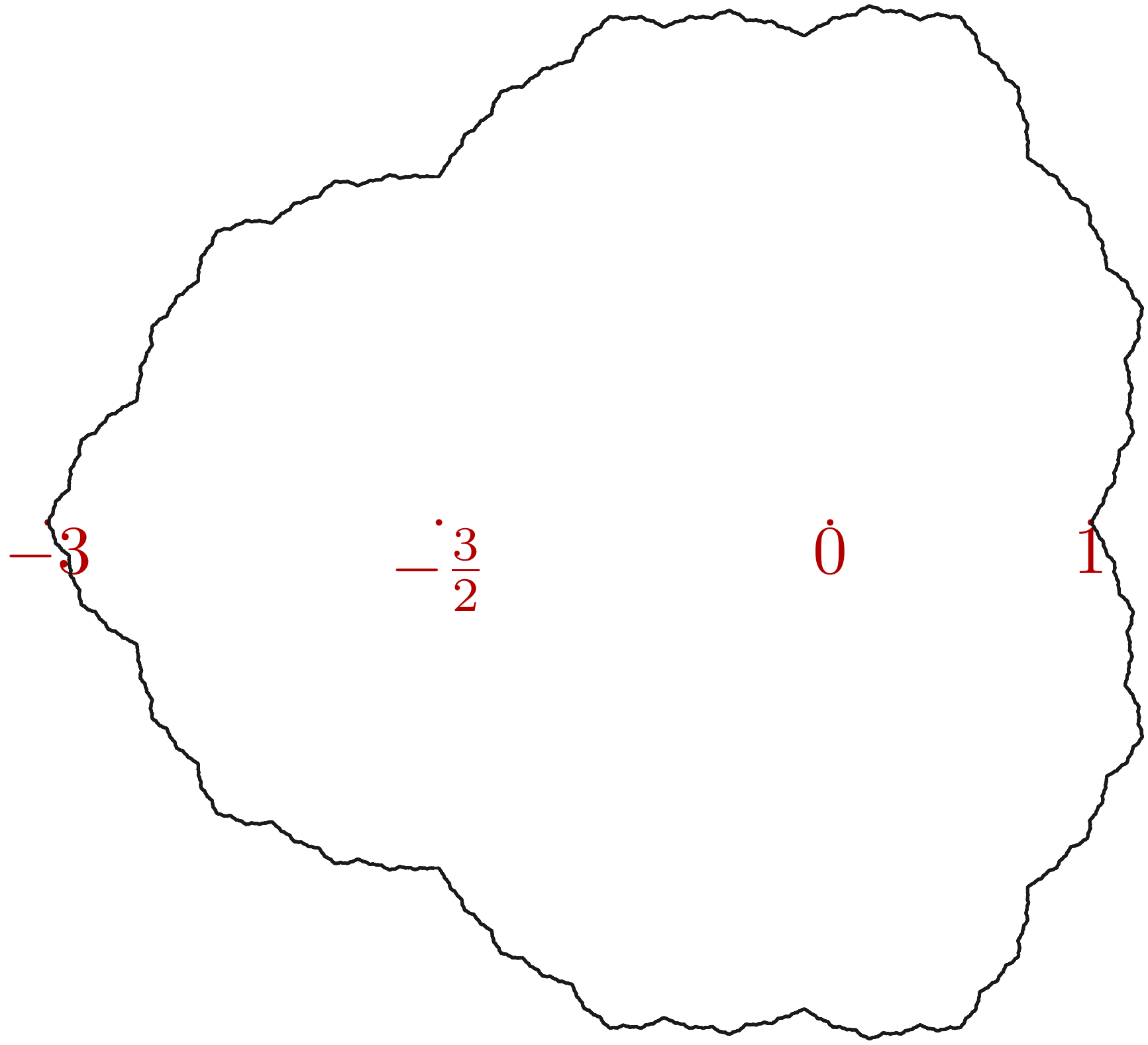}}
	\caption{Julia set for the polynomial $P(z)=\frac14(z+2z^2+z^3)$. It is the boundary of the open, connected, and bounded set, the filled Julia set, which is the domain of definition for the analytic function $\Phi$, see \er{004} and \er{005}. Note that for other polynomials, Julia sets can be disconnected but always bounded (or empty).  }\lb{fig1}
\end{figure}
For the numerical computation of $\Phi(z)$ one may rewrite \er{004} as a recurrent relation
\[\lb{006}
 \Phi_0(z)=z,\ \ \ \Phi_{n+1}(z)=\Phi_{n}(z)\lt(1+\frac{2\Phi_{n}(z)}{4^n}+\frac{\Phi_{n}(z)^2}{16^n}\rt),\ \ n\ge0.
\]
Then $\Phi_n(z)\to\Phi(z)$ exponentially fast, and uniformly on any fixed compact subset of the filled Julia set. Substituting Tailor expansion
\[\lb{007}
 \Phi(z)=\vp_1z+\vp_2z^2+\vp_3z^3+...\ \ \ (\vp_1=1)
\]
along with $P(z)=\frac{z}4(z+1)^2$ into \er{005} we obtain the recurrent formula for $\vp_i$:
\[\lb{008}
 \vp_1=1,\ \ \ \vp_{i+1}=\frac{4}{1-4^{-i}}\sum_{k=0}^{i-1}4^{k-i}\binom{2(i-k)}{k+1}\vp_{i-k},\ \ i\ge1, 
\]
where $\binom{a}{b}$ are binomial coefficients. We assume $\binom{a}{b}=0$ if $b>a$. For numerical implementations, it is useful to apply $\exp\circ\ln$ with some straightforward simplifications to the coefficients before $\vp_{i-k}$ in \er{008}. Combining results of Theorem \ref{T1} along with \er{004}, \er{007}, and \er{008} we obtain the following Corollary. 

\begin{corollary}\lb{C1}
For the densities $\r_i$ of the cells with weigh $i$ in the exponential ordered (OSC) and disordered (DSC) stochastic compression models, the following identities hold
\[\lb{009}
{\rm OSC}:\ \ \ \frac{\r_i}{\r_1}=i,\ \ N>\log_2i;\ \ \ \ \ {\rm DSC}:\ \ \ \frac{\r_i}{\r_1}\to\vp_i,\ \ N\to\iy,
\]
where $\vp_i$ are given by \er{008}, and $N$ is the number of compression steps.
\end{corollary}

It is important to note that $\r_1=4^{-N}$ for both: ordered and disordered cases. First few values $\vp_i$ are given in Tab. \ref{Tab2}. A comparison of the values from Tabs. \ref{Tab2} and \ref{Tab1}.(d) illustrates the results of Corollary \ref{C1}.

\begin{table}[h]
\begin{center}
    \begin{tabular}{ | l | l | l | l | l | l | l | l | l |}
    \hline
    i & $1$ & $2$ & $3$ & $4$  & $5$ & $6$ & $7$ & $8$ \\ \hline
    $\vp_i$ & $1.00000$  & $2.66667$  & $3.91111$ & $5.55344$  & $7.05507$  & $8.26885$  & $9.86538$ & $11.41518$ \\ 
    \hline
    \end{tabular}
 \end{center}
\caption{Tailor coefficients of $\Phi(z)$ computed by \er{008}.}\lb{Tab2}
\end{table}

Now, we will focus on the most interesting part of the article - approximation of $\vp_i$. At first, consider the branch $P^{-1}(z)$ which satisfies $P^{-1}(1)=1$. For the cubic polynomial $P(z)=\frac{z+2z^2+z^3}4$, the inverse function can be computed explicitly
\[\lb{010}
 P^{-1}(z)=\frac{\sqrt[3]{54z+1+\sqrt{(54z+1)^2-1}}+\sqrt[3]{54z+1-\sqrt{(54z+1)^2-1}}-2}3.
\]
The mapping $P^{-1}$ is an analytic injection on, at least, the complex plane with the slit $\C\sm(-\iy,0]$, since
\[\lb{011}
 P^{-1}=A_5\circ A_4\circ A_3\circ A_2\circ A_1,
\]
where
\[\lb{012}
 A_1(z)=54z+1,\ \ A_2(z)=z+\sqrt{z^2-1},\ \ A_3(z)=\sqrt[3]{z},\ \ A_4(z)=z+\frac1{z},\ \ A_5(z)=\frac{z-2}3
\]
are known conformal mappings.
Now, we define the mapping $\Psi$ by analogy with $\Phi$, see \er{004},
\[\lb{013}
 \Psi(z)=\lim_{N\to\iy}2^N(\underbrace{P^{-1}\circ...\circ P^{-1}}_{N}(z)-1).
\]
The definition is correct, since $1$ is an attracting point with the multiplier $\frac12$: $P^{-1}(1)=1$ and $(P^{-1})'(1)=\frac12$. Function $\Psi$ is an analytic injection on $\C\sm(-\iy,0]$ as a limit composition of analytically injective and contractive mappings $P^{-1}$. It satisfies the functional equation
\[\lb{014}
 \Psi(P(z))=2\Psi(z),\ \ \ \Psi'(1)=1,
\]
see \er{013}. Using
\[\lb{015}
 P(z)-1=(z-1)\frac{4+3z+z^2}{4}
\]
along with \er{013}, we obtain the recurrent sequence of functions
\[\lb{016}
 \Psi_0(z)=z-1,\ \ \ \Psi_{n}(z)=\Psi_{n-1}(z)\frac{8}{4+3\underbrace{P^{-1}\circ...\circ P^{-1}}_{n}(z)+\underbrace{P^{-1}\circ...\circ P^{-1}}_{n}(z)^2},
\]
which converges $\Psi_n(z)\to\Psi(z)$ uniformly on any compact subset of $\C\sm(-\iy,0]$. The convergence is exponentially fast, since $\underbrace{P^{-1}\circ...\circ P^{-1}}_{n}(z)\to1$ exponentially fast with the factor $\frac12$, recall that $P^{-1}(1+z)=1+\frac{z}2+O(z^2)$. Formula \er{016} is especially good for numerical computations. Denote
\[\lb{017}
 \Theta(z)=\Phi(z)\Psi(z)^2.
\]
Function $\Theta(z)$ is analytic on the intersection of the domains of definition of functions $\Phi(z)$ and $\Psi(z)$, i.e. on $\mJ_P\sm[-3,0]$, where $\mJ_{P}$ is the Julia set for $P(z)=\frac{z+2z^2+z^3}4$, see Fig. \ref{fig1}. Using \er{005} and \er{014} we conclude that $\Theta(z)$ satisfies the functional equation
\[\lb{018}
 \Theta(P(z))=\Theta(z).
\]
Function $\Theta(z)$ is well defined on the open interval $(0,1)$ which is a subset of its domain of definition $\mJ_{P}\sm[-3,0]$, see above. Denote
\[\lb{019}
 \vt_{\max}=\max_{z\in(0,1)}\Theta(z),\ \ \ \vt_{\min}=\min_{z\in(0,1)}\Theta(z).
\]
Due to \er{018}, for the evaluation of $\vt_{\max}$ and $\vt_{\min}$ it is enough to compute maximum and minimum of $\Theta(z)$ on the interval $[P(a),a]$ with any $a\in(0,1)$. We take $a=\frac12$ and the corresponding interval $[\frac9{32},\frac12]$, which is far enough from the singularities of functions $\Phi$, $\Psi$, and, hence, $\Theta$. Accurate computations based on \er{006} and \er{016} show the wonderful result
\[\lb{020}
 \vt_{\max}=1.46491046...,\ \ \ \vt_{\min}=1.46491015...,
\]
i.e. $\Theta(z)$ is ``almost constant" for $z\in(0,1)$. However, $\Theta(z)$ oscillates between $\vt_{\min}$ and $\vt_{\max}$ infinitely many times for $z\in(0,1)$. In the final chapter of {\bf Remark 2} at the end of the article, we provide the detailed explanation of the fact why $\Theta(z)$ is ``almost constant". Denoting
\[\lb{021}
 \vt_0=\frac{\vt_{\max}+\vt_{\min}}2=1.464910...,
\]
we obtain the best constant approximation of $\Theta(z)$ in the uniform norm
\[\lb{022}
 \Theta(z)=\vt_0+\ve(z)=1.464910...+\ve(z),\ \ \ |\ve(z)|<10^{-6},\ \ z\in(0,1),
\]
where $\ve(z)$ is analytic on the same domain as $\Theta(z)$. Using \er{017} and \er{022}, we get
\[\lb{023}
 \Phi(z)=\frac{\vt_0}{\Psi(z)^2}+\frac{\ve(z)}{\Psi(z)^2}.
\]
Differentiating \er{014} twice at $z=1$, we obtain the Tailor series $\Psi(z)=(z-1)-\frac{5(z-1)^2}8+O((z-1)^3)$, which leads to
\[\lb{024}
 \frac1{\Psi(z)^2}=\frac1{(z-1)^2}+\frac{5}{4(z-1)}+O(1),\ \ z\to1.
\]
Combining  \er{022}-\er{024}, we obtain
\[\lb{025}
 \Phi(z)=\frac{\vt_0}{(1-z)^2}-\frac{5\vt_0}{4(1-z)}+{\rm "small\ terms"}.
\]
Taking Tailor expansions of both sides in \er{025}, we obtain
\[\lb{026}
 \vp_i\approx(i+1)\vt_0-\frac54\vt_0=(i-\frac14)\vt_0.
\]
Approximation \er{026} is really good, see Fig. \ref{fig2}, but, unfortunately, we can not state that $\vp_i-(i-\frac14)\vt_0$ is $o(i)$. More accurate approximations require more difficult reasoning that goes beyond our current needs. 
This is just an entrance to the rabbit hole leading to the fractal abyss. 
The main statements announced in the Abstract follows from \er{021}, \er{026}, and Corollary \ref{C1}. Very briefly, further analysis can be based on the Stirling asymptotic of factorials in binomial coefficients in \er{008}, in the same way as in the known first variants of the proof of CLT. The main impact in RHS of \er{008} is for $k\in[i/2-c\sqrt{i},i/2+c\sqrt{i}]$ with $c>0$. Roughly speaking, this fact implies $\vp_{i}\approx 2\vp_{i/2}$, which, in turn, leads to $\vp_i\approx i\t(\log_2i)$ with $1$-periodic function $\t$. This $1$-periodic function is a near-constancy oscillation $\t\approx\vt_0$. Perhaps, some ideas from \cite{O} can be helpful in this analysis, but I am not sure about that. Some relevant results are given in Theorem \ref{T3} below.
\begin{figure}[h]
    \centering
    \begin{subfigure}[b]{0.45\textwidth}
        \includegraphics[width=\textwidth]{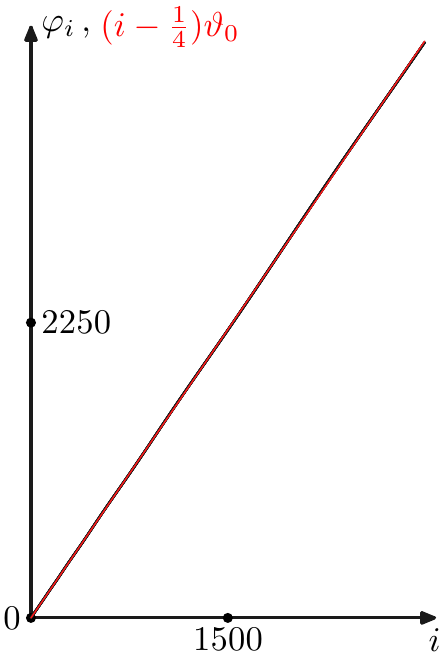}
        \caption{}
    \end{subfigure}
    \begin{subfigure}[b]{0.45\textwidth}
        \includegraphics[width=\textwidth]{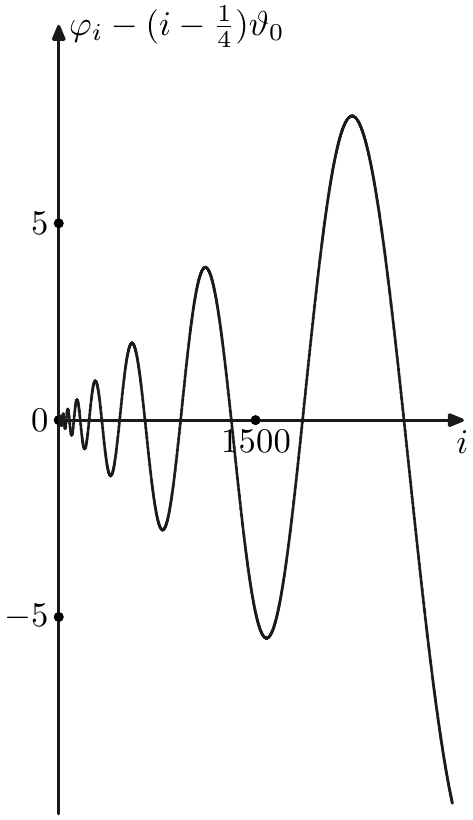}
        \caption{}
    \end{subfigure}
    \caption{Tailor coefficients $\vp_i$ (black curve) and their approximation $(i-\frac14)\vt_0$ (red curve) are plotted in (a), the difference $\vp_i-(i-\frac14)\vt_0$ is plotted in (b). For the computation of $\vp_i$ formula \er{008} is used.}\lb{fig2}
\end{figure}

Let us provide an extension of Theorem \ref{T1} and Corollary \ref{C1} related to the inhomogeneous chain $\Z$. We consider the simple inhomogeneous case where cells of the chain at the initial state can have weight $1$ or $0$ (empty cells). Suppose that the initial density of non-empty cells is $1\ge p>0$. In other words, two-valued independent random variables $\zeta_j$, such that $\zeta_j=1$ with the probability $p$ and $\zeta_j=0$ with the probability $1-p$, are placed at each cell $j\in\Z$.

\begin{theorem}\lb{T2} Suppose that the chain $\Z$ has an initial simple inhomogeneous state with the density $1\ge p>0$ of non-empty cells. Then densities $\r_i$ of cells with weight $i$ after $N$ steps of compression have the form
\[\lb{027}
 \r_i=F^{(i)}(1-p)\frac{p^i}{i!},\ \ \ i\ge0,
\]
where
\[\lb{028}
 F(z)=\frac{z(z^{2^{N}}-1)^2}{4^N(z-1)^2}\ \ \ or\ \ \ F(z)=\underbrace{P\circ...\circ P}_{N}(z)
\]
for the ordered or disordered stochastic compression models respectively. The limits of relative densities are
\[\lb{029}
 \lim_{N\to\iy}\frac{\r_i}{\r_1}=\frac{i+1-p}{2-p},\ \ \ \lim_{N\to\iy}\frac{\r_i}{\r_1}=\frac{\Phi^{(i)}(1-p)\,p^{i-1}}{\Phi'(1-p)\,i!},\ \ \ i\ge0
\] 
in the ordered and disordered cases respectively.
\end{theorem}

Let us illustrate results, namely \er{029}, with numerical examples having the same setup as in Tab. \ref{Tab1}, but with the inhomogeneous initial chain, where half of the cells are empty, and half of the cells have weight $1$. The empty and non-empty cells are distributed randomly in the initial chain. The numerical results are in good agreement with the theory, see Tab. \ref{Tab3}. The derivatives of $\Phi$ in \er{029} are computed by taking derivatives of the Tailor series, see \er{007} and \er{008}.  Note that at each step of compression non-empty cells do not disappear. The linear growth of densities in the ordered case is obvious. The analysis of growth in the disordered case is more complex. It can be performed in the same way as for \er{026}. However, even simple differentiation of \er{025}, already gives a very good approximation
\[\lb{030}
 \frac{\Phi^{(i)}(1-p)\,p^{i-1}}{\Phi'(1-p)\,i!}\approx\frac{(i+1-\frac54p)\vt_0}{p^3\Phi'(1-p)},\ \ \ i\ge0.
\]
For example, for $p=0.5$ and $i=0,...,50$ the maximal deviation RHS from LHS in \er{030} is less than $0.063$  by absolute value, while RHS increase linearly as $0.729...(i+\frac38)$.

\begin{table}[h]
\begin{center}
    \begin{tabular}{ | l | l | l | l | l | l | l | l | l |}
    \hline
    i & $0$ & $1$ & $2$ & $3$  & $4$ & $5$ & $6$ & $7$ \\ \hline
    $N_i$ & $5006$  & $15112$  & $25259$ & $34504$  & $44530$  & $55101$  & $65140$ & $74244$ \\ \hline
    $\frac{N_i}{N_1}$ & $0.33126$ & $1.00000$ & $1.67145$ & $2.28322$  & $2.94666$  & $3.64618$  & $4.31048$ & $4.91292$ \\  \hline
    $\frac{i+1-p}{2-p}$ & $0.33333$ & $1.00000$ & $1.66667$ & $2.33333$  & $3.00000$  & $3.66667$  & $4.33333$ & $5.00000$ \\
    \hline
    \end{tabular}
\\
\bigskip
(a) OSC, $8$ steps
\bigskip
\\
    \begin{tabular}{ | l | l | l | l | l | l | l | l | l |}
    \hline
    i & $0$ & $1$ & $2$ & $3$  & $4$ & $5$ & $6$ & $7$ \\ \hline
    $N_i$ & $6250$  & $20088$  & $34538$ & $49073$  & $63641$  & $78082$  & $93337$ & $106675$ \\ \hline
    $\frac{N_i}{N_1}$ & $0.31113$ & $1.00000$ & $1.71933$ & $2.44290$  & $3.16811$  & $3.88700$  & $4.64641$ & $5.31038$ \\ 
     \hline
    $\frac{\Phi^{(i)}(1-p)\,p^{i-1}}{\Phi'(1-p)\,i!}$ & $0.31495$ & $1.00000$ & $1.73321$ & $2.46170$  & $3.19152$  & $3.92065$  & $4.65000$ & $5.38002$ \\
    \hline
    \end{tabular}
\\    
    \bigskip
    (b) DSC, $8$ steps 
\end{center}
\caption{The quantities $N_i$ of cells with weights $i$ in a segment of inhomogeneous chain of length $2^{26}\cdot 5^4\approx42\cdot10^9$ after $8$ steps of ordered and disordered stochastic compression, in numerical simulations. The initial density of cells with weight $1$ in the inhomogeneous chain is $p=\frac12$. The ratio $N_i/N_1$ is compared with theoretical values \er{029}.}\lb{Tab3}
\end{table}

Above, we discuss mostly the densities of cells with relatively small weights. In view of \cite{K1}, this is the most interesting case. At the same time, it is possibly to look on $\r_i$ globally in the following sense. Formula \er{002} suggests a proper rescaling. Define
\[\lb{031}
 \r(x)=\lim_{N\to\iy}2^N\r_{[x2^N]},\ \ \ x\in[0,\iy),
\]
where $\r_i$ are the densities of cells $i$ after $N$ steps of compression, and $[a]$ denotes an integer part of the real number $a$. For the ordered case, \er{002} leads to
\[\lb{032}
 \r(x)=\ca x,& x\in[0,1],\\
           2-x, & x\in[1,2],\\
           0,& x\ge2.
       \ac    
\]
In Fig. \ref{fig3}, rescaled densities \er{031} are compared in ordered and disordered cases. It is seen how the presence of random permutations in the disordered case can change the distribution of densities. Note that both integrals $\int_0^{+\iy}\r(x)dx=1$. The notable difference appears near $x=0$ - the random permutations increases the ``linear" growth of densities. Above, we already discussed the similar increase $146$ percent.

\begin{figure}[ht]
	\center{\includegraphics[width=0.9\linewidth]{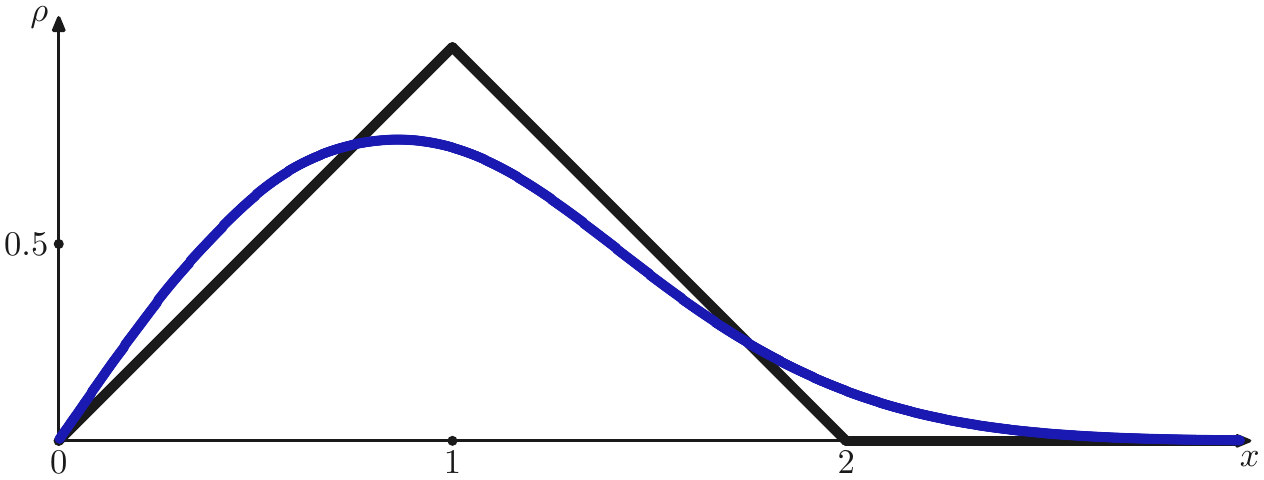}}
	\caption{Rescaled densities $\r(x)$, defined by \er{031}, in the ordered (black curve) and disordered (blue curve) cases are plotted. For the computation of points of blue curve, $9$ iterations of $P$ are used.}\lb{fig3}
\end{figure}

We denote $\mathbf{i}=\sqrt{-1}$. For the disordered case, it is possible to compute $\r(x)$ implicitly, by using its characteristic function $\chi(t)$ with the help of \er{003}, Riemann sums for integrals, and the asymptotic $e^x=1+x+o(x)$, $x\to0$:
\begin{multline}\lb{033}
 \chi(t)=\int_{-\iy}^{\iy}e^{\mathbf{i}tx}\r(x)dx=\lim_{N\to\iy}2^{-N}\sum_{i=1}^{\iy}2^N\r_ie^{\frac{\mathbf{i}ti}{2^N}}=\lim_{N\to\iy}\underbrace{P\circ...\circ P}_{N}(e^{\frac{\mathbf{i}t}{2^N}})=\\
 \lim_{N\to\iy}\underbrace{P\circ...\circ P}_{N}(1+\frac{\mathbf{i}t}{2^N})=\Pi(\mathbf{i}t),\ \ \ t\in\C,
\end{multline}
where $\Pi$ is an entire function satisfying a Poincar\'e-type functional equation
\[\lb{034}
 P(\Pi(z))=\Pi(2z),\ \ \Pi'(0)=1,\ \ \ z\in\C,
\] 
see the last identity in \er{033}. The existence of $\Pi$ easily follows from the fact that $1$ is a repelling point for the polynomial $P(z)=\frac 14(z+2z^2+z^3)$, namely $P(1+z)=1+2z+o(z)$, $z\to0$. A good computational procedure for evaluation of $\Pi(z)$ is given in \er{313}.  Applying the inverse Fourier transform to $\chi(t)$, we can recover $\r(x)$ in the disordered case. We do not give the details of the calculations, but the recovered density 
\[\lb{invF}
 \r(x)=\frac1{2\pi}\int_{-\iy}^{\iy}e^{-\mathbf{i}tx}\Pi(\mathbf{i}t)dt
\] 
is indistinguishable from the blue curve presented in Fig. \ref{fig3}. The alternative formula expressed $\r(x)$ through a lacunary-type series is given in \er{321}. Differentiating \er{034} at $z=0$ and using $\Pi(0)=P(1)=1$, it is easy to write explicit values $\Pi^{(n)}(0)$, which, in turn, coincide with the corresponding moments of $\r$, since the moments are derivatives of $\chi(-\mathbf{i}t)$ at $t=0$, see \er{033}. Thus, we have
\[\lb{035}
 \int_0^{\iy}x\r(x)dx=\Pi'(0)=1,\ \ \int_0^{\iy}x^2\r(x)dx=\Pi''(0)=\frac54,\ \ \int_0^{\iy}x^3\r(x)dx=\Pi'''(0)=\frac{87}{48}
\]
and so on. Again, all the quantities \er{035} are in good agreement with numerical computations. In the end, one may ask: we define various Schr\"oder and Poincar\'e-type functions. Are there any other connections between them? Yes, using \er{014} and \er{034}, or direct definitions \er{013} and \er{033}, it is not difficult to see that
\[\lb{036}
 \P(\Pi(z))=\Pi(\P(z))=z,
\]
at least in the neighborhood of $z=1$. These functions are mutually inverse to each other. Of course, an analytic continuation of \er{036} to sufficiently large regions exists. Let us return to our main course: the asymptotic of $\r(x)$ when $x\to0$. It is possible to prove the next theorem.
\begin{theorem}\lb{T3}
The following asymptotic holds
\[\lb{res}
 \r(x)=x\s(\log_2x)+o(x),\ \ \ x\to0,
\]
where $\s$ is $1$-periodic continuous function defined by the formula
\[\lb{res1}
 \s(x)=2^{-x}\r(2^x)+\frac{2^{-x}}{4\pi}\int_{-\iy}^{+\iy}H(-\mathbf{i}t2^{x+1})(2\Pi(\mathbf{i}t)^2+\Pi(\mathbf{i}t)^3)dt,\ \ \ x\in\R,
\]
where the entire function $H$ satisfies
\[\lb{res2}
 H(z)=\sum_{n=2}^{\iy}\frac{z^n}{n!(2^{n-1}-1)},\ \ \ 
 H(2z)=2H(z)+2(e^z-1-z),\ \ \ z\in\C.
\]
The convergence of the integral in \er{res1} is ensured by the facts that $\Pi(\mathbf{i}t)=O(|t|^{-2})$ and $H(\mathbf{i}t)=O(|t|\ln|t|)$ for large by modulus $t\in\R$. Moreover, there are a couple of identities for the average of $\s(x)$:
\begin{multline}\lb{C007}
 \int_{0}^1\s(x)dx=\frac1{2\ln2}\int_0^1\Psi(z)^2(4z+3z^2)dz=\frac{-1}{\ln2}\int_{-\iy}^{0}t(2\Pi(t)^2+\Pi(t)^3)dt=\\ 
 \frac{-1}{\ln2}\int_{0}^{+\iy}t\Re(2\Pi(\mathbf{i}t)^2+\Pi(\mathbf{i}t)^3)dt=\frac{-2}{\pi\ln2}\int_{0}^{+\iy}(\ln t)t\Im(2\Pi(\mathbf{i}t)^2+\Pi(\mathbf{i}t)^3)dt.
\end{multline}
The convergence of the third integral in \er{C007} is due to $\Pi(t)=O(t^{-2})$ for large $t$ with $\Re t\le0$.
\end{theorem}
Using ideas explained in \er{317}-\er{319}, one can show that the integral in \er{res1} is differentiable by $x$. This means that the remainder $o(x)$ in \er{res} is smooth for $x\ne0$. Computations show that
$\s(x)\approx\int_0^1\s(x)dx\approx1.464910...$, see Fig. \ref{figS}, as it is also expected from the results provided at the beginning of this article. A little bit informal, but detailed analysis of the connection between these constants is done in Remark 2 at the end of this article, see \er{C009} and \er{C010}. For numerical approaches, it is useful to take into account that the best rate of convergence is shown by using the last formula in \er{C007}. 
\begin{figure}[ht]
	\center{\includegraphics[width=0.7\linewidth]{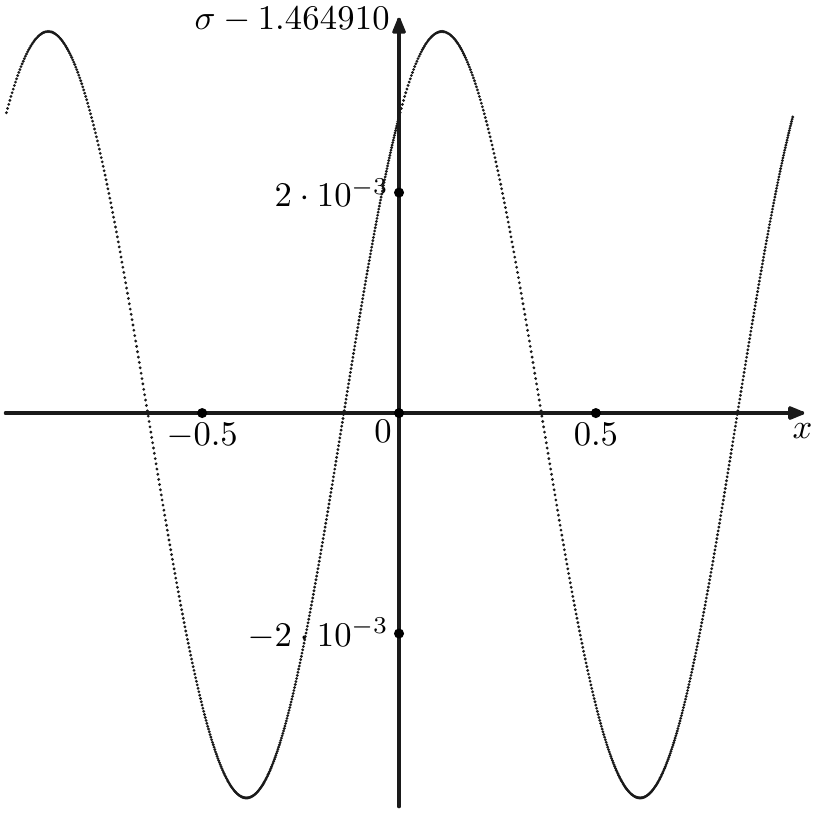}}
	\caption{The $1$-periodic function $\s(x)-1.464910$, see \er{res}, computed by using \er{res1}, \er{res2}, and \er{invF}.}\lb{figS}
\end{figure}
 
For the inhomogeneous chain with the initial density $p$ of non-empty cells, \er{027} leads to
\[\lb{037}
 \sum_{i=0}^{+\iy}\r_iz^i=F(1-p+pz),
\]  
where $F$ is given by \er{028}. Thus, in the disordered case, for the rescaled densities $\r_p(x)$, defined similarly to \er{031}, the corresponding characteristic function 
\[\lb{038}
 \chi_p(t)=\chi(pt),
\]
as it is not difficult to verify with the help of \er{037} and \er{033}. Hence, in the disordered case, we have  
\[\lb{039}
 \r_p(x)=\frac1p\r\lt(\frac xp\rt)
\] 
by the known property of the Fourier transform, see \er{038}. Using \er{037}, first identity in \er{028}, and the characteristic function approach, it is easy to check that \er{039} is also true for the ordered case.

{\bf Remark 1.} For the disordered models, where cells in the chain merge left and right with different probabilities, say $r$ and $1-r$, most of the results mentioned in this Section remains true if we replace $P(z)=\frac{z+2z^2+z^3}4$ with $P_r(z)=r(1-r)z+(r^2+(1-r)^2)z^2+r(1-r)z^3$. Of course, the Julia set for $P_r(z)$ depends on $r$. The results related to one-dimensional disordered compression models also carry over to multidimensional cases, with the only difference that the polynomial $P$, see \er{003}, can have a degree greater than $3$. At the same time, the analysis of the ordered case with $r\not=\frac12$ and in a multidimensional case seems more complex.

Now we are going to the proof of Theorems \ref{T1} and \ref{T2}.

{\section{Proof of Theorems \ref{T1} and \ref{T2}}\lb{sec2}}

\subsection{Proof of Theorem \ref{T1}. Ordered compression.}

Suppose that the system lives $N$ steps. Consider a segment of the length $2n+1$. Define the polynomial
\[\lb{100}
 P_N(x_0,...,x_{2n})=\sum_{(r_0,...,r_{2n})\in\N^{2n+1}}\r_N(r_0,...,r_{2n})x_0^{r_0}...x_{2n}^{r_{2n}},
\]
where $\r_N(r_0,...,r_{2n})$ is the density of the sequence $(r_0,...,r_{2n})$ in the entire model after $N$ steps. Note that we omit the second parameter $2n+1$ writing simply $P_N$, since the number of variables indicates this parameter uniquely. The first goal is to find recurrence relations for
the polynomials $P_N$. Introduce 
\[\lb{101}
 \mM_{2n}=\lt\{(y_0,...,y_{4n+2}):\ \ca y_0=\es\ {\rm or}\ 0,\\
                                                y_{2j+1}=j,\ {\rm for}\ j=0,...,2n,\\
                                                y_{2j}={j-1}\ {\rm or}\ {j},\ {\rm for}\                                                          j=1,...,2n,\\
                                                y_{4n+2}=\es\ {\rm or}\ {2n} \ac\rt\}.
\]
In other words, the numbers $j$ are placed in odd cells of vectors from $\cM$. The numbers in even cells coincide with one of the number located in the neighboring left or right cell. In addition, the first and the last elements of vectors in $\cM$ can be equal to the empty element $\es$. The number of vectors in this set is $\#\cM_{2n}=2^{2n+2}$. Further, we need also the number of vectors in some of its subsets
\begin{multline}\lb{102}
 \#\{{\bf y}\in\cM_{2n}:\ y_0=\es,\ y_{4n+2}=\es\}=
 \#\{{\bf y}\in\cM_{2n}:\ y_0=0,\ y_{4n+2}=\es\}=\\
 \#\{{\bf y}\in\cM_{2n}:\ y_0=\es,\ y_{4n+2}=2n\}=\#\{{\bf y}\in\cM_{2n}:\ y_0=0,\ y_{4n+2}=2n\}=2^{2n}.
\end{multline}
The structure of the vectors from $\cM({\bf x})$ reflects the dynamics of our system. Namely, the polynomials $P_N$, $N\ge1$ satisfy
\[\lb{103}
 P_N(x_0,...,x_{2n})=2^{-2n-2}\sum_{{\bf m}=(m_j)\in\mM_{2n}}P_{N-1}(x_{m_0},x_{m_1},...,x_{m_{4n+2}}),
\]
where $x_{\es}$ is assumed to be equal $1$. RHS of \er{103} describes all the possible sequences that transforms to ${\bf x}$ in one step. For example,
\[\lb{104}
 P_N(x_0)=\frac14(P_{N-1}(1,x_0,1)+P_{N-1}(x_0,x_0,1)+P_{N-1}(1,x_0,x_0)+P_{N-1}(x_0,x_0,x_0))=...,
\]
which can be continued further up to linear combinations of $P_0$. We fix $N$. Let us associate the polynomial $Q_i(x)$ which counts the number of $x_0$ in formulas \er{104} multiplied by the corresponding coefficients also appearing in \er{104}, i.e.
\[\lb{105}
 Q_N(x)=x,\ \ \ Q_{N-1}(x)=\frac14(x+2x^2+x^3)
\]
and so on. Using \er{101} along with \er{103}, it is seen that the arguments of $P_i$ in \er{104} consist of three connected segments $(1,...,1,x_0,...,x_0,1,...,1)$ maybe of zero length some of them. Using this fact and \er{102}, we deduce that
\[\lb{106}
 Q_{i-1}(x)=\cL Q_i(x),\ \ i\ge1,
\]
where the linear operator $\cL$ acts on the basic polynomials by the following rule
\[\lb{107}
 \cL x^j=\frac14(x^{2j-1}+2x^{2j}+x^{2j+1}),\ \ j\ge1.
\]
By induction, it is easy to check that
\[\lb{108}
 \cL^ix=\frac1{4^i}\sum_{j=1}^{2^i}jx^j+\frac1{4^i}\sum_{j=2^i+1}^{2^{i+1}-1}(2^{i+1}-j)x^j,
\]
since
\begin{multline}\lb{ind1}
 \cL\lt(\frac1{4^i}\sum_{j=1}^{2^i}jx^j+\frac1{4^i}\sum_{j=2^i+1}^{2^{i+1}-1}(2^{i+1}-j)x^j\rt)=\frac1{4^{i+1}}\sum_{j=1}^{2^i}j(x^{2j-1}+
 2x^{2j}+x^{2j+1})+\\
 \frac1{4^{i+1}}\sum_{j=2^i+1}^{2^{i+1}-1}(2^{i+1}-j)(x^{2j-1}+
 2x^{2j}+x^{2j+1})=\frac1{4^{i+1}}\sum_{j=1}^{2^{i+1}}jx^j+\frac1{4^{i+1}}\sum_{j=2^{i+1}+1}^{2^{i+2}-1}(2^{i+2}-j)x^j.
\end{multline}
Thus, using \er{105}, \er{106}, and \er{108}, we obtain
\[\lb{109}
 Q_0(x)=\frac1{4^N}\sum_{j=1}^{2^N}jx^j+\frac1{4^N}\sum_{j=2^N+1}^{2^{N+1}-1}(2^{N+1}-j)x^j.
\]
Remembering that $Q_0$ counts entries of $x_0$ in the polynomials $P_0$ with the proper coefficients appearing in \er{104}, and using the fact that, by definition, $P_0(x_0,...,x_{2n})=\prod_{j=0}^{2n}x_j$ for any $n$, we deduce that
\[\lb{110}
 P_N(x_0)=Q_0(x_0).
\]
Identities \er{109}, \er{110} and the definition of the polynomial $P_N$ finishes the proof.

\subsection{Proof of Theorem \ref{T1}. Disordered compression.}

Introduce the polynomials
\[\lb{111}
 P_N(z)=\sum_{n=1}^{\iy}\r_{n,N}z^n,
\]
where $\r_{n,N}$ is the density of the cell with weight $n$ in the chain after $N$ steps of compression. The probability of merging three elements in one cell is $1/4$. The probability of merging only two elements in one cell is $1/2$, since it may happens in two cases: only left or only right element moves into the central cell. The probability of the case when the central element leaves alone is again $1/4$. Thus, we obtain the following recurrent relation
\begin{multline}\lb{112}
 P_{N+1}(z)=\sum_{n=1}^{\iy}\lt(\frac14\r_{n,N}+\frac12\sum_{n_1+n_2=n}\r_{n_1,N}\r_{n_2,N}+\frac14\sum_{n_1+n_2+n_3=n}\r_{n_1,N}\r_{n_2,N}\r_{n_3,N}\rt)z^n=\\
 \frac14P_{N}(z)+\frac12P_N(z)^2+\frac14P_N(z)^3=P(P_N(z))=\underbrace{P\circ...\circ P}_{N+1}(z),
\end{multline} 
since $P_0(z)$ is equal to $z$.

\subsection{Proof of Theorem \ref{T2}.} Let $F(z)$ be a generating function for the homogeneous compression models after $N$ steps
\[\lb{113}
 F(z)=\sum_{i\ge 1}\r_i^{\rm h}z^i,
\]
where densities $\r_i^{\rm h}$ of the cells with weight $i$ are given by \er{002} in ordered or \er{003} in disordered cases. Then, in the inhomogeneous case, after $N$ steps of compression the densities $\r_i$ of the cells with weight $i$ will be
\[\lb{114}
 \r_j=\sum_{i\ge 1}\r_i^{\rm h}\mathbb{P}\{\zeta_{a_{1i}}+...+\zeta_{a_{ii}}=j\}=\sum_{i\ge 1}\frac{\r_i^{\rm h}i!p^j(1-p)^{i-j}}{j!(i-j)!}=F^{(j)}(1-p)\frac{p^j}{j!},
\]  
since $a_{xy}$ are all different and the corresponding $\zeta_{a_{xy}}$ are all independent two-valued random variables. Thus, \er{027} and second identity in \er{028} follows from \er{113} and \er{003}.  

Now, note that
\[\lb{115}
 S(n,z):=\sum_{i=1}^niz^i=z\lt(\sum_{i=1}^nz^i\rt)'=z\lt(\frac{z^{n+1}-z}{z-1}\rt)'=\frac{nz^{n+2}-(n+1)z^{n+1}+z}{(z-1)^2}
\]
Let us compute explicitly \er{113} for $\r_i^{\rm h}$ given by \er{002}. In this case, we have
\begin{multline}\lb{116}
 F(z)=4^{-N}\sum_{i=1}^{2^N}iz^i+4^{-N}\sum_{i=2^N+1}^{2^{N+1}-1}(2^{N+1}-i)z^i=S(2^N,z)+z^{2^{N+1}}S(2^N-1,z^{-1})=\\
 \frac{2^Nz^{2^N+2}-(2^N+1)z^{2^N+1}+z}{4^N(z-1)^2}+z^{2^{N+1}}\frac{(2^N-1)z^{-2^N-1}-2^Nz^{-2^N}+z^{-1}}{4^N(z^{-1}-1)^2}=\\
 \frac{2^Nz^{2^N+2}-(2^N+1)z^{2^N+1}+z}{4^N(z-1)^2}+\frac{(2^N-1)z^{2^N+1}-2^Nz^{2^N+2}+z^{2^{N+1}+1}}{4^N(z-1)^2}=\\
 \frac{z^{2^{N+1}+1}-2z^{2^N+1}+z}{4^N(z-1)^2}=\frac{z(z^{2^{N}}-1)^2}{4^N(z-1)^2},
\end{multline}
where \er{115} is used. First identity in \er{028} follows from \er{116}. First identity in \er{116} leads to
\[\lb{117}
 4^NF(z)\to\sum_{i=1}^{\iy}iz^i=\frac1{(1-z)^2}-\frac1{1-z},\ \ \ N\to\iy.
\]
The convergence is uniform on any compact subset of the unit ball. Thus, the derivatives also converge on such subsets uniformly
\[\lb{118}
 4^NF^{(j)}(z)\to\frac{(j+1)!}{(1-z)^{2+j}}-\frac{j!}{(1-z)^{1+j}},\ \ \ N\to\iy.
\]
Identities \er{027}, \er{117}, and \er{118} lead to the first identity in \er{029}, since the factor $4^N$ disappears in the ratio. The same arguments applied to \er{004} lead to the second identity in \er{029}. 

\section{Proof of Theorem \ref{T3}. \lb{sec3}}

We have already shown that $\Theta(z)$, see \er{017}, is bounded and separated from $0$ for $z\in(0,1)$. Using the same arguments based on \er{018}, it is not difficult to show that $\Theta(z)$ is bounded and and separated from $0$ in a neighborhood of $z=0$. Thus, using \er{017} along with $\Psi(z)=z+o(z)$, $z\to0$, we obtain that $\Psi(z)$ is of the order $1/\sqrt{z}$ up to a bounded and separated from $0$ multiplier, when $z\to0$. Using this fact, it is seen that the inverse function $\Pi=\Psi^{-1}$ satisfies $|\Pi(\mathbf{i}t)|\le a/t^2$ with constant $a\ge0$ for large positive and negative $t\in\R$. Thus, there is no absolute convergence of $\int_{-\iy}^{\iy} t\Pi(\mathbf{i}t)dt$, and, hence, we cannot state that $\r'(0)$ exists, see \er{invF}. Let us define another integral
\[\lb{300}
 B:=\frac1{2\pi \mathbf{i}}\int_{-\iy}^{\iy}t(2\Pi(\mathbf{i}t)^2+\Pi(\mathbf{i}t)^3)dt,
\] 
which converges absolutely, since $|\Pi(\mathbf{i}t)|\le a/t^2$ for large by modulus $t\in\R$. Now, let us change the variables $t\leftrightarrow 2t$ in \er{invF} and use \er{034}:
\begin{multline}\lb{301}
\r(x)=\frac1{2\pi}\int_{-\iy}^{\iy}e^{-\mathbf{i}t2x}\Pi(\mathbf{i}t)dt=\frac1{\pi}\int_{-\iy}^{\iy}e^{-\mathbf{i}2tx}\Pi(\mathbf{i}2t)dt=\\
\frac1{\pi}\int_{-\iy}^{\iy}e^{-\mathbf{i}t2x}\frac{\Pi(\mathbf{i}t)+2\Pi(\mathbf{i}t)^2+\Pi(\mathbf{i}t)^3}4dt=\frac{\r(2x)}{2}+R(x),
\end{multline}
where
\[\lb{302}
 R(x)=\frac1{4\pi}\int_{-\iy}^{+\iy}e^{-\mathbf{i}t2x}(2\Pi(\mathbf{i}t)^2+\Pi(\mathbf{i}t)^3)dt.
\]
Due to $|\Pi(\mathbf{i}t)|\le a/t^2$, the function $R(x)$ has at least two derivatives. The derivative of $R$ at $x=0$ coincides with \er{300}, and, hence $R'(0)=B$. The value $R(0)$ is $0$, since $\r(0)=0$. Using $R(0)=0$ and $R'(0)=B$, we can write
$R(x)=Bx+O(x^2)$ for $x\to0$, which with \er{301} gives
\[\lb{303}
 \lt|\frac{\r(x)}{x}+B\log_2x-\lt(\frac{\r(2x)}{2x}+B\log_22x\rt)\rt|\le C|x|,\ \ \ 0<x\le2
\] 
for some $C\ge0$. For $n\ge0$, let us define continuous functions
\[\lb{304}
 r_n(x)=\frac{\r(\frac x{2^n})}{\frac x{2^n}}+B\log_2\frac{x}{2^n},\ \ \ x\in[1,2].
\]
The functions are continuous, since $\r$ is continuous as a Fourier transform of rapidly decaying function $\Pi(\mathbf{i}t)$. Using \er{303} and \er{304}, we obtain
\[\lb{305}
 \lt|r_n(x)-r_{n-1}(x)\rt|\le C2^{-n+1},\ \ \ x\in[1,2],\ \ \ n\ge1,
\]
which, in turn, leads to
\[\lb{306}
 \lt|r_m(x)-r_{n}(x)\rt|\le C\sum_{j=n}^{m-1}2^{-j}\le C2^{-n+1},\ \ \ x\in[1,2],\ \ \ m>n\ge1.
\]
Estimates \er{306} show that $\{r_n(x)\}$ is a Cauchy sequence of continuous functions. It converges in the uniform norm to a continuous function, i.e.
\[\lb{307}
 r_n(x)\rightrightarrows r(x)
\] 
uniformly on the interval $x\in[1,2]$, where $r(x)$ is some continuous function such that $r(1)=r(2)$, since $r_{n}(2)=r_{n-1}(1)$. Consider continuous $\s(x)=r(2^x)$, $x\in[0,1]$. Since $\s(0)=\s(1)$, we can extend $\s(x)$ to $1$-periodic continuous function defined for all $x\in\R$. Using \er{304}, \er{307} and $1$-periodicity of $\s$, we obtain
\begin{multline}\lb{308}
 \r(\frac x{2^n})=\frac{x}{2^n}r(x)-B\frac{x}{2^n}\log_2\frac{x}{2^n}+o(\frac{x}{2^n})=\frac{x}{2^n}\s(\log_2x)-B\frac{x}{2^n}\log_2\frac{x}{2^n}+o(\frac{x}{2^n})=\\
 \frac{x}{2^n}\s(\log_2x-n)-B\frac{x}{2^n}\log_2\frac{x}{2^n}+o(\frac{x}{2^n})=\frac{x}{2^n}\s(\log_2\frac{x}{2^n})-B\frac{x}{2^n}\log_2\frac{x}{2^n}+o(\frac{x}{2^n}),
\end{multline}
which, after the change of variables $t=\frac{x}{2^n}$, gives 
\[\lb{309}
 \r(t)=t\s(\log_2t)-Bt\log_2t+o(t),\ \ \ t\to0.
\]
Everything is good, but $B=0$. Indeed, for convenience, we change the variable $z=\Pi(\mathbf{i}t)$ or $t=-\mathbf{i}\Psi(z)$, see \er{036}. Then \er{300} becomes
\[\lb{310}
 B=\frac{-1}{2\pi\mathbf{i}}\oint_{\g}\Psi(z)(2z^2+z^3)\Psi'(z)dz=\frac{1}{4\pi\mathbf{i}}\oint_{\g}\Psi(z)^2(4z+3z^2)dz,
\]
where the contour $\g=\Pi(\mathbf{i}\R)$ is depicted in Fig. \ref{fig4}. The contour $\g$ is continuous and connected, since analytic function $\Pi(it)\to0$ for real $t\to\pm\iy$. We only need to check that $\oint_{|z|=\ve}\Psi(z)^2(4z+3z^2)dz\to0$ for $\ve\to0$. It is enough to show that the integrand in \er{310} is bounded. It is true, since
\[\lb{311}
 \Psi(z)^2(4z+3z^2)=\Theta(z)(4+3z)\frac{z}{\Phi(z)},
\]  
see \er{017}. We already discussed that $\Theta$ is bounded near $z=0$, and $\frac{z}{\Phi(z)}\to1$ for $z\to0$. We also ensure the correctness of integration by parts in \er{310}, since $|2z^2+z^3|\le|4z+3z^2|$ near $z=0$. The analytic function $\Psi$ has no singularities inside the domain bounded by the contour $\g$. Thus $B=0$ in \er{310} by the Cauchy's residue theorem, and we obtain the main result \er{res} from \er{309}.
It is useful to note that $\Psi$ conformally maps the domain bounded by the contour $\gamma$ onto the left half-plane $\{t:\ \Re t<0\}$, since $\Psi((0,1))=(-\iy,0)$ and $\Psi(\g)=\mathbf{i}\R$. In particular, this means that $\Pi(t)=O(|t|^{-2})$ for large $t$ with $\Re t\le0$, since $\Pi=\Psi^{-1}$ and $\Psi(z)=O(|z|^{-\frac12})$ near $z=0$.
\begin{figure}[ht]
	\center{\includegraphics[width=0.5\linewidth]{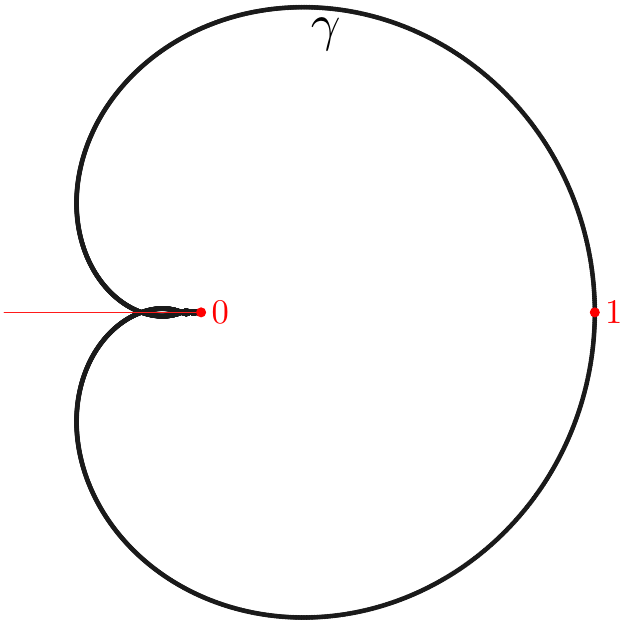}}
	\caption{The contour $\g=\Pi(\mathbf{i}\R)$ in the complex plane. The domain of definition of $\Psi$ is $\C\sm(-\iy,0]$. The corresponding slit is depicted by red.}\lb{fig4}
\end{figure}
Finally, let us discuss some explicit expressions for $\s(x)$. Formula \er{301} and \er{302} along with $R(0)=R'(0)=0$ gives
\begin{multline}\lb{314}
 \frac{\r(\frac{x}{2^n})}{\frac{x}{2^n}}-\frac{\r(x)}{x}=\sum_{j=1}^n\frac{R(\frac{x}{2^j})}{\frac{x}{2^j}}=\sum_{j=1}^n\frac{R(\frac{x}{2^j})-R(0)-R'(0)\frac{x}{2^j}}{\frac{x}{2^j}}=\\
 \sum_{j=1}^n\frac1{4\pi}\int_{-\iy}^{+\iy}\frac{e^{-\mathbf{i}t2\frac{x}{2^j}}-1+\mathbf{i}t2\frac{x}{2^j}}{\frac{x}{2^j}}(2\Pi(\mathbf{i}t)^2+\Pi(\mathbf{i}t)^3)dt\xrightarrow[n\to\iy]{}\\
 \frac1{4\pi x}\int_{-\iy}^{+\iy}H(-\mathbf{i}t2x)(2\Pi(\mathbf{i}t)^2+\Pi(\mathbf{i}t)^3)dt,
\end{multline}
where
\[\lb{315}
 H(z)=\sum_{j=1}^{\iy}2^j(e^{\frac{z}{2^j}}-1-\frac{z}{2^j})=\sum_{j=1}^{\iy}2^j\sum_{n=2}^{\iy}\frac{z^n}{2^{jn}n!}=\sum_{n=2}^{\iy}\frac{z^n}{n!(2^{n-1}-1)}
\]
is an entire function of the first exponential order of growth $|H(z)|\le4 e^{\frac{|z|}2}$, satisfying functional equation
\[\lb{316}
 H(2z)=2H(z)+2(e^z-1-z),
\]
which follows from the first identity in \er{315}. To justify the change $\int$ ans $\sum$ and the convergence in \er{314}, the exponential growth of $F$ is not enough, since $\Pi(it)=O(t^{-2})$ only. We need more accurate estimates based on \er{316}, namely
\[\lb{317}
 \frac{H(2z)}{2z}=\frac{H(z)}{z}+\frac{e^z-1-z}{z},
\]
which leads to
\[\lb{318}
 \lt|\frac{H(2\mathbf{i}t)}{2t}\rt|\le\lt|\frac{H(\mathbf{i}t)}{t}\rt|+2,\ \ \ t\in\R\sm(-2,2).
\]
Taking any $C\ge2$ such that $|\frac{H(\mathbf{i}t)}{t}|\le C\log_2|t|$ for $|t|\in[2,4]$, we deduce from \er{318} by induction the estimate
\[\lb{319}
 \lt|\frac{H(\mathbf{i}t)}{t}\rt|\le C\log_2|t|
\]
for any $t$ such that $|t|\ge2$. The estimate \er{319} is already enough for the convergence and all the manipulations in \er{314}. Combining \er{res} with \er{314}, we obtain
\[\lb{320}
 \s(\log_2 x)=\frac{\r(x)}{x}+\frac1{4\pi x}\int_{-\iy}^{+\iy}H(-\mathbf{i}t2x)(2\Pi(\mathbf{i}t)^2+\Pi(\mathbf{i}t)^3)dt,
\]
which with \er{314}-\er{316}, and \er{319} give \er{res1} and \er{res2}.

Now, we fix some large $R>0$. We have
\begin{multline}\lb{C001}
 A_1(R):=\int_0^1\frac{2^{-x}}{2\pi}\int_{-R}^Re^{-\mathbf{i}t2^x}\Pi(\mathbf{i}t)dtdx=\int_0^1\frac{2^{-x}}{2\pi}\int_{-R}^R\sum_{n=0}^{+\iy}\frac{(-\mathbf{i}t2^x)^n}{n!}\Pi(\mathbf{i}t)dtdx=\\
 \frac1{4\pi\ln2}\int_{-R}^R\Pi(\mathbf{i}t)dt-\frac1{2\pi}\int_{-R}^R\mathbf{i}t\Pi(\mathbf{i}t)dt+\sum_{n=2}^{+\iy}\frac1{2\pi\ln2}\int_{-R}^R\frac{(2^{n-1}-1)(-\mathbf{i}t)^n}{(n-1)n!}\Pi(\mathbf{i}t)dt.
\end{multline}
Using \er{res1}, \er{res2}, and \er{034}, where $P(z)=\frac{z+2z^2+z^3}4$, we have
\begin{multline}\lb{C002}
 A_2(R):=\int_0^1\frac{2^{-x}}{4\pi}\int_{-R}^{+R}H(-\mathbf{i}t2^{x+1})(2\Pi(\mathbf{i}t)^2+\Pi(\mathbf{i}t)^3)dtdx=\\
 \int_0^1\frac{2^{-x}}{4\pi}\int_{-R}^{+R}H(-\mathbf{i}t2^{x+1})(4\Pi(2\mathbf{i}t)-\Pi(\mathbf{i}t))dtdx=\\
 \int_0^1\frac{2^{-x}}{4\pi}\int_{-R}^R\sum_{n=2}^{+\iy}\frac{(-\mathbf{i}t2^{x+1})^n}{n!(2^{n-1}-1)}(4\Pi(2\mathbf{i}t)-\Pi(\mathbf{i}t))dtdx=\\
 \sum_{n=2}^{+\iy}\frac1{4\pi\ln2}\int_{-R}^R\frac{(-\mathbf{i}t)^n2^n}{(n-1)n!}(4\Pi(2\mathbf{i}t)-\Pi(\mathbf{i}t))dt.
\end{multline}
Taking the sum of \er{C001} and \er{C002}, we obtain
\begin{multline}\lb{C003}
 A_1(R)+A_2(R)=
 \frac1{4\pi\ln2}\int_{-R}^R\Pi(\mathbf{i}t)dt-\frac1{2\pi}\int_{-R}^R\mathbf{i}t\Pi(\mathbf{i}t)dt+\\
 -\sum_{n=2}^{+\iy}\frac1{2\pi\ln2}\int_{-R}^R\frac{(-\mathbf{i}t)^n}{(n-1)n!}\Pi(\mathbf{i}t)dt+\sum_{n=2}^{+\iy}\frac1{\pi\ln2}\int_{-R}^R\frac{(-\mathbf{i}t)^n2^n}{(n-1)n!}\Pi(2\mathbf{i}t)dt=\\
 \frac1{4\pi\ln2}\int_{-R}^R\Pi(\mathbf{i}t)dt-\frac1{2\pi}\int_{-R}^R\mathbf{i}t\Pi(\mathbf{i}t)dt+\frac1{2\pi\ln2}(\int_{-2R}^{2R}-\int_{-R}^{R})\sum_{n=2}^{+\iy}\frac{(-\mathbf{i}t)^n}{(n-1)n!}\Pi(\mathbf{i}t)dt=\\
 \frac1{4\pi\ln2}\int_{-R}^R\Pi(\mathbf{i}t)dt-\frac1{2\pi}\int_{-R}^R\mathbf{i}t\Pi(\mathbf{i}t)dt+\frac1{2\pi\ln2}(\int_{-2R}^{2R}-\int_{-R}^{R})\int_0^{it}\frac{\mathbf{i}t(e^{\zeta}-1+\zeta)d\zeta}{\zeta^2}\Pi(\mathbf{i}t)dt=\\
 \frac1{4\pi\ln2}\int_{-R}^R\Pi(\mathbf{i}t)dt-\frac1{2\pi}\int_{-R}^R\mathbf{i}t\Pi(\mathbf{i}t)dt+\frac1{2\pi\ln2}(\int_{-2R}^{-R}+\int_{R}^{2R})\mathbf{i}t(\ln\mathbf{i}t+D+\ve(t))\Pi(\mathbf{i}t)dt=\\
 \frac1{4\pi\ln2}\int_{-R}^R\Pi(\mathbf{i}t)dt+\frac{D}{2\pi\ln2}(\int_{-2R}^{2R}-\int_{-R}^{R})\mathbf{i}t\Pi(\mathbf{i}t)dt+\frac1{2\pi\ln2}(\int_{-2R}^{-R}+\int_{R}^{2R})\mathbf{i}t\ve(t)\Pi(\mathbf{i}t)dt\\
 \frac1{2\pi\ln2}(\int_{-2R}^{2R}-\int_{-R}^{R})\mathbf{i}t(\ln\mathbf{i}t)\Pi(\mathbf{i}t)dt-\frac1{2\pi}\int_{-R}^R\mathbf{i}t\Pi(\mathbf{i}t)dt=\\
 \frac1{4\pi\ln2}\int_{-R}^R\Pi(\mathbf{i}t)dt+\frac{D}{2\pi\ln2}\int_{-R}^{R}\mathbf{i}t(2\Pi(\mathbf{i}t)^2+\Pi(\mathbf{i}t)^3)dt+\frac1{2\pi\ln2}(\int_{-2R}^{-R}+\int_{R}^{2R})t\ve(t)\Pi(\mathbf{i}t)dt+\\
 \frac{1}{2\pi\ln2}\int_{-R}^{R}\mathbf{i}t(\ln\mathbf{i}t)(2\Pi(\mathbf{i}t)^2+\Pi(\mathbf{i}t)^3)dt\xrightarrow[R\to+\iy]{}\frac{1}{2\pi\ln2}\int_{\R}\mathbf{i}t(\ln\mathbf{i}t)(2\Pi(\mathbf{i}t)^2+\Pi(\mathbf{i}t)^3)dt.
\end{multline}
In this long derivation \er{C003} we use a lot of facts, such that \er{034}, \er{300} with $B=0$, $\r(0)=\frac1{2\pi}\int_{\R}\Pi(\mathbf{i}t)dt=0$, $\Pi(\mathbf{i}t)=O(t^{-2})$ where real $t\to\pm\iy$, and
\[\lb{C004}
 \sum_{n=2}^{+\iy}\frac{(-\mathbf{i}t)^n}{(n-1)n!}=\int_0^{it}\frac{\mathbf{i}t(e^{\zeta}-1+\zeta)d\zeta}{\zeta^2}=\mathbf{i}t(\ln\mathbf{i}t+D+\ve(t)),\ \ \ \ve(t)=O(t^{-1})
\]
for real $t\to\pm\iy$ and constant $D$. Now, using \er{C003}, and \er{C001} with \er{invF}, and \er{C002} with \er{res1}, we deduce that
\begin{multline}\lb{C005}
 \int_{0}^1\s(x)dx=\frac{1}{2\pi\ln2}\int_{\R}\mathbf{i}t(\ln\mathbf{i}t)(2\Pi(\mathbf{i}t)^2+\Pi(\mathbf{i}t)^3)dt=\\
 \frac{1}{2\pi\ln2}\int_{\R}\mathbf{i}t(\ln|t|+\frac{\mathbf{i}\pi\sign t}2)(2\Pi(\mathbf{i}t)^2+\Pi(\mathbf{i}t)^3)dt.
\end{multline}
If we take another branch of logarithm $\ln_1\mathbf{i}t=\ln|t|-\frac{\mathbf{i}\pi\sign t}2$ with the slit $[0,+\iy)$ then $\ln_1\Psi(z)$ will be analytic inside the domain bounded by the contour $\gamma$, see Fig. \ref{fig4}, since $\Psi((0,1])=(-\iy,0]$. Thus, using the same arguments as in the calculation of $B$, see \er{300} and \er{310}, we obtain
\begin{multline}\lb{C006}
 \frac{1}{2\pi\ln2}\int_{\R}\mathbf{i}t(\ln|t|-\frac{\mathbf{i}\pi\sign t}2)(2\Pi(\mathbf{i}t)^2+\Pi(\mathbf{i}t)^3)dt=
 \frac{1}{2\pi\ln2}\int_{\R}\mathbf{i}t(\ln_1\mathbf{i}t)(2\Pi(\mathbf{i}t)^2+\Pi(\mathbf{i}t)^3)dt=\\
 \frac{1}{2\pi\mathbf{i}\ln2}\oint_{\g}\Psi(z)(\ln_1\Psi(z))(2z^2+z^3)\Psi'(z)dz=\\
 \frac{-1}{2\pi\mathbf{i}\ln2}\oint_{\g}\frac{\Psi(z)^2(2\ln_1\Psi(z)-1)}4(4z+3z^2)dz=0.
\end{multline}
Again, integration by parts in \er{C006}, and the absence of difficulties near $z=0$ is ensured by the same methods as in the calculation of $B$, see \er{300} and \er{310}, namely $\Psi(z)=O(|z|^{-\frac12})$ near $z=0$, and $\Psi(z)\ln_1\Psi(z)\to0$ near $z=1$. Combining \er{C005} and \er{C006}, remembering $\Pi(-\mathbf{i}t)=\ol{\Pi(\mathbf{i}t)}$, 
we obtain the corresponding equations in \er{C007}. To obtain other equations, we take the standard logarithm $\ln$ with the slit $(-\iy,0]$ in \er{C006}, so that it coincides with \er{C005}, and it is not $0$, since $\Psi((0,1))$ is exactly the slit for $\ln$, and the interval $(0,1)$ to which we would like to shrink $\g$ lies inside the domain bounded by $\gamma$ now:
\begin{multline}\lb{C008}
 \int_{0}^1\s(x)dx=\frac{1}{2\pi\ln2}\int_{\R}\mathbf{i}t(\ln\mathbf{i}t)(2\Pi(\mathbf{i}t)^2+\Pi(\mathbf{i}t)^3)dt=\\
 \frac{-1}{2\pi\mathbf{i}\ln2}\oint_{\g}\frac{\Psi(z)^2(2\ln\Psi(z)-1)}4(4z+3z^2)dz=\\\frac{1}{2\pi\mathbf{i}\ln2}(\int_{1+\mathbf{i}0}^{0+\mathbf{i}0}+\int_{0-\mathbf{i}0}^{1-\mathbf{i}0})\frac{\Psi(z)^22\ln\Psi(z)}4(4z+3z^2)dz=\frac1{2\ln2}\int_0^1\Psi(z)^2(4z+3z^2)dz=\\
 =\frac1{2\ln2}\int_{-\iy}^{0}t^2(4\Pi(t)+3\Pi(t)^2)\Pi'(t)dt=\frac{-1}{\ln2}\int_{-\iy}^{0}t(2\Pi(t)^2+\Pi(t)^3)dt,
\end{multline}
which give other equations in \er{C007}. In \er{C008} we use the fact that $\Pi(t)=O(|t|^{-2})$ for $\Re t\le0$, which follows from $C_1|z|^{-\frac12}\le|\Psi(z)|\le C_2|z|^{-\frac12}$ for some constants $C_1,C_2>0$, in the same way as already discussed above.

{\bf Remark 2.} Computations show that $\int_0^1\s(x)dx=1.464910...$. Are there any connections with the similar value related to the function $\Theta(z)$, see \er{017}-\er{020}. Yes, it is. Using \er{C008} and varying $0<a<1$ and $n\in\N$, we have
\begin{multline}\lb{C009}
\int_0^1\s(x)dx=\frac1{\ln4}\int_0^1\Psi(z)^2(4z+3z^2)dz=\frac1{\ln4}\int_{\underbrace{P\circ...\circ P}_{n}(a)}^a\Psi(z)^2(4z+3z^2)dz+O_1(a,n)=\\
\frac1{\ln4}\int_{\underbrace{P\circ...\circ P}_{n}(a)}^a\Psi(z)^2(4P'(z)-1)dz+O_1(a,n)=\frac1{\ln4}\int_{\underbrace{P\circ...\circ P}_{n+1}(a)}^{P(a)}\Psi(z)^2dz+\\
-\frac1{\ln4}\int_{\underbrace{P\circ...\circ P}_{n}(a)}^a\Psi(z)^2dz+O_1(a,n)=\frac1{\ln4}\int_{\underbrace{P\circ...\circ P}_{n+1}(a)}^{\underbrace{P\circ...\circ P}_{n}(a)}\Psi(z)^2dz+\\
-\frac1{\ln4}\int_{P(a)}^{a}\Psi(z)^2dz+O_1(a,n)=
\frac1{\ln4}\int_{\underbrace{P\circ...\circ P}_{n+1}(a)}^{\underbrace{P\circ...\circ P}_{n}(a)}\Psi(z)^2\Phi'(z)dz+O_2(a,n)+\\-\frac1{2\ln2}\int_{P(a)}^{a}\Psi(z)^2dz+
O_1(a,n)=\frac1{\ln4}\int_{4^{-n-1}\Phi(a)}^{4^{-n}\Phi(a)}\Psi(\Phi^{-1}(z))^2dz+O_2(a,n)+\\-\frac1{2\ln2}\int_{P(a)}^{a}\Psi(z)^2dz+
O_1(a,n)=\int_{\log_4\Phi(a)-n-1}^{\log_4\Phi(a)-n}4^z\Psi(\Phi^{-1}(4^z))^2dz+O_2(a,n)+\\-\frac1{2\ln2}\int_{P(a)}^{a}\Psi(z)^2dz+
O_1(a,n)=\int_{\log_4\Phi(a)-n-1}^{\log_4\Phi(a)-n}\Theta(\Phi^{-1}(4^z))dz+O_2(a,n)+
\\-\frac1{2\ln2}\int_{P(a)}^{a}\Psi(z)^2dz+
O_1(a,n),
\end{multline} 
and
\[\lb{C010}
 \lim_{a\to1}\lim_{n\to\iy}\sum_{j=1}^2O_j(a,n)=0,\ \ \ \lim_{a\to1}\frac1{2\ln2}\int_{P(a)}^{a}\Psi(z)^2dz=0,
\]
where we use \er{004}, \er{005}, $\Phi'(0)=P(1)=1$, $\underbrace{P\circ...\circ P}_{n}(a)\to0$ for $n\to\iy$, and \er{014} with $\Psi(z)=O(|z|^{-\frac12})$. Using the already mentioned functional equations for $\Phi$ and $\Theta$, it is easy to check that $\Theta(\Phi^{-1}(4^z))$ is $1$-periodic analytic function defined on some neighborhood of $\R$. Identities \er{C009} and \er{C010} show that the averages of $\s(x)$ and $\Theta(\Phi^{-1}(4^z))$ over their periods are identical. However, there is a still interesting problem about connections between $L^{\iy}$- and $L^{1}$-norms of $\Theta(\Phi^{-1}(4^z))$ over period. It seems that they are ``practically" identical and equal to $1.464910...$.

A direct computation of the periodic function $\Theta(\Phi^{-1}(4^z))$ depends on the accurate computation of $\Phi^{-1}$. We have already the accurate procedure for the computation of $\Theta(z)$, see \er{017}, \er{018} along with \er{016} and  \er{006}. Let us discuss the computation of
\[\lb{C011}
 \Phi^{-1}(z)=\lim_{N\to\iy}\underbrace{P^{-1}\circ...\circ P^{-1}}_{N}(4^{-N}z),
\]
see \er{004}. It is very inconvenient to calculate \er{C011}, since $4^{-N}z$ is very small argument, which should be computed very precisely. Using
\[\lb{C012}
 P^{-1}(\frac z4)=\lt(\frac{3}{T(z)+T(z)^{-1}+1}\rt)^2z,\ \ \ T(z)=(54z+1+\sqrt{(54z+1)^2-1})^{\frac13},
\]
see \er{010}, we introduce
\[\lb{C013}
 \wt\Phi_N(z):=\underbrace{P^{-1}\circ...\circ P^{-1}}_{N}(4^{-N}z)=
 \wt\Phi_{N-1}\lt(\lt(\frac{3}{T(4^{-N}z)+T(4^{-N}z)^{-1}+1}\rt)^2z\rt).
\]
Using the second identity in \er{C013}, we obtain a rapid procedure for computing the limit $\Phi^{-1}(z)=\lim_{N\to\iy}\wt\Phi_N(z)$, see \er{C011}, since the factor before $z$ in the last term of \er{C013} tends to $1$ exponentially fast for $N\to\iy$. This fact allows us to avoid an accurate computation of $4^{-N}z$ for large $N$ and moderate $z$. The second identity in \er{C013} is recurrent, that is a little bit different from other fast procedures mentioned above, see \er{016} and \er{006}. We also use this idea below to obtain a fast algorithm for the computation of $\Pi(z)$, see \er{313}.

To finalize this remark, let us discuss the domain of definition of the $1$-periodic function $\Theta(\Phi^{-1}(4^z))$. As the function $\Psi$, see the explanation in \er{010}-\er{017}, the function $\Theta(z)$ is analytic in some neighborhood of $z=0$ without the cut along the negative values of the real axis. At the same time $\Phi^{-1}(z)=z+o(z)$ is analytic near $z=0$. Thus, taking $z=-R+\mathbf{i}\vp$ with $R,\vp>0$, we see that $\Phi^{-1}(4^z)$ reaches the cut along the negative values of the real axis when $\vp\to\frac{\pi}{\ln 4}$ for $R\to\iy$. Hence, due to the $1$-periodicity and the real-axis symmetry, the function 
$\Theta(\Phi^{-1}(4^z))$ is analytic in the strip $\{z:\ |{\rm Im z}|<\frac{\pi}{\ln 4}\}$. Applying a simple transformation 
$
 K(z):=\Theta(\Phi^{-1}(4^{\frac{\ln z}{2\pi\mathbf{i}}})) 
$,
we obtain the function analytic in the ring $\{z:\ e^{-\frac{\pi^2}{\ln 2}}<|z|<e^{\frac{\pi^2}{\ln 2}}\}$. Moreover, the computations show that $|K(z)|<2$ in this ring, see Fig. \ref{fig5}. Hence, the Cauchy estimates yield that the $n$-th Laurent coefficient of $K(z)$ is less than $2e^{-|n|\frac{\pi^2}{\ln 2}}$ by modulus. Taking into account the fact that $e^{\frac{\pi^2}{\ln 2}}>10^6$, we may see why $K(z)$ is ``almost" constant that leads to the fact that $\Theta(x)$ is also ``almost" constant for $x\in(0,1)$.

\begin{figure}[ht]
	\center{\includegraphics[width=0.7\linewidth]{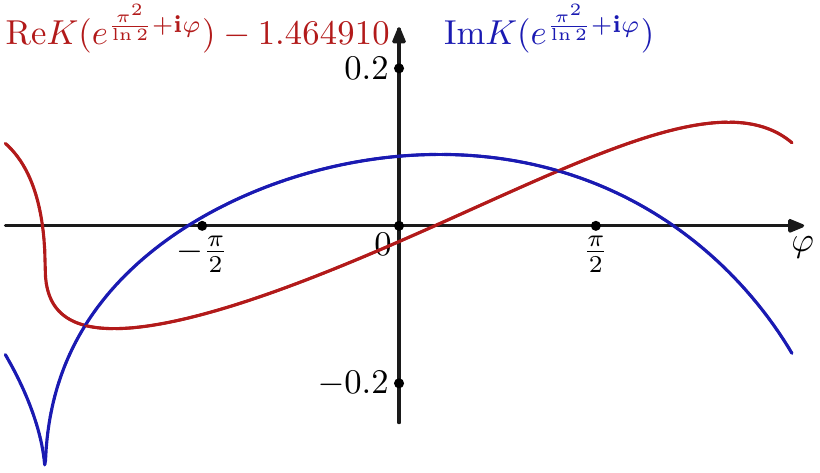}}
	\caption{The real (red) and imaginary (blue) part of $K(z)$ at the circle $|z|=e^{\frac{\pi^2}{\ln 2}}$. The plot of $K(z)$ on another edge circle $|z|=e^{\frac{\pi^2}{\ln 2}}$ is similar. Some singular behavior at one point is visible.}\lb{fig5}
\end{figure}

{\bf Remark 3.} Let us discuss some alternative to \er{invF} formulas for the computation of $\r(x)$. Again, using \er{301}, \er{302} and the fact that $\r(x)\to0$ for $x\to\iy$ as a Fourier transform of a rapidly decaying function, we obtain
\[\lb{321}
 \r(x)=\frac1{2\pi}\int_{-\iy}^{\iy}G(-\mathbf{i}t2x)(2\Pi(\mathbf{i}t)^2+\Pi(\mathbf{i}t)^3)dt,
\]
where
\[\lb{322}
 G(z)=\sum_{j=0}^{\iy}\frac{e^{2^jz}}{2^j}
\]
is an analytic in the left half-plane $\{z:\ \Re z<0\}$ function, bounded $|G(z)|\le2$ there, and continuous up to the boundary $\{z:\ \Re z=0\}$, satisfying the functional equation
\[\lb{323}
 G(2z)=2G(z)-2e^z.
\]
In fact, $\r(x)x^n\to0$, $x\to\iy$ for any $n\ge0$, since its (inverse) Fourier transform $\Pi^{(n)}(\mathbf{i}t)\to0$ at least as $O(|t|^{-2})$ for real $t$ as can be shown by using the functional equation \er{034} in a similar way as in \er{317}-\er{319} with the help of already proven fact $\Pi(\mathbf{i}t)=O(|t|^{-2})$ for real $t$.

{\bf Remark 4.} In this final remark, let us discuss a numerical procedure that allows us to compute $\Pi(z)$, see \er{033} and \er{034}. At first, note that 
\[\lb{312}
 P(1+\frac{z}{2^N})=\frac{1+\frac{z}{2^N}+2(1+\frac{z}{2^N})^2+(1+\frac{z}{2^N})^3}4=1+\frac{z(1+\frac{5z}{2^{N+3}}+\frac{z^2}{2^{2N+3}})}{2^{N-1}}.
\] 
Thus, using \er{312}, we can rewrite the composition limit from \er{034} by introducing the following functions
\[\lb{313}
 \Pi_0(z)=1+z,\ \ \ \Pi_n(z)=\Pi_{n-1}\lt(z+\frac{5z^2}{2^{n+3}}+\frac{z^3}{2^{2n+3}}\rt),\ \ \ n\ge1.
\]
Then $\Pi_n(z)\to\Pi(z)$ uniformly  and exponentially fast on any compact set. The type of a composition formula \er{313} differs from that one of \er{006} and \er{016}.

\section*{Acknowledgements} 
This paper is a contribution to the project M3 of the Collaborative Research Centre TRR 181 "Energy Transfer in Atmosphere and Ocean" funded by the Deutsche Forschungsgemeinschaft (DFG, German Research Foundation) - Projektnummer 274762653. 


\begin{thebibliography}{9}
\bibitem{K1}
A. A. Kutsenko,  ``Dynamics of group formation". {\it https://arxiv.org/abs/2112.12733}, 2021.

\bibitem{F} W. Feller, ``An introduction to probability theory and its applications", Vol. 1, 3nd Ed., John Willey \& Sons, 1970.

\bibitem{M} J. Milnor, ``Dynamics in one complex variable", Princeton NJ, Princeton Univ. Press, 2006.

\bibitem{K2}
A. Kutsenko,  ``An entire function connected with the approximation of the golden ratio". {\it Am. Math. Monthly}, {\bf 127}, 820-826, 2020.

\bibitem{CG}
O. Costin and G. Giacomin,  ``Oscillatory critical amplitudes in hierarchical models and the Harris function of branching processes". {\it J. Stat. Phys.}, {\bf 150}, 471–486, 2013.

\bibitem{H}
T. E. Harris, ``Branching processes", {\it Ann. Math. Statist.},  {\bf 41}, 474-494, 1948.

\bibitem{BB} 
J. D. Biggings and N. H. Bingham, ``Near-constancy phenomena in branching processes". {\it Math. Proc. Cambridge Phil. Soc.}, {\bf 110}, 545-558, 1991.

\bibitem{DIL}
B. Derrida, C. Itzykson, and J. M. Luck,  ``Oscillatory critical amplitudes in hierarchical models". {\it Comm. Math. Phys.}, {\bf 94}, 115–132, 1984.

\bibitem{O}
A. M. Odlyzko,  ``Periodic oscillations of coefficients of
power series that satisfy functional equations". {\it Adv. Math.}, {\bf 44}, 180–205, 1982.

\end{thebibliography}
\end{document}